\documentclass{article}

\usepackage[margin=1.2in]{geometry}

\usepackage{amsmath,amssymb,amsfonts,graphicx,amsthm,nicefrac,mathtools,bm}
\usepackage{hyperref}
\usepackage[numbers]{natbib}
\usepackage{bbm}
\usepackage[ruled,vlined,linesnumbered,resetcount]{algorithm2e}
\usepackage{multirow}
\usepackage{graphicx}
\graphicspath{ {images/} }
\newtheorem{theorem}{Theorem}
\usepackage{rotating}
\usepackage{hhline}
\usepackage[affil-it]{authblk}

\newcommand{\ymax}{y_{\max}}
\newcommand{\R}{\mathbb{R}}
\newcommand{\muh}{\widehat{\mu}}
\newcommand{\rh}{\widehat{r}}
\newcommand{\muf}{\muh_{\textnormal{fast}}}
\newcommand{\mus}{\muh_{\textnormal{slow}}}
\newcommand{\rf}{\rh_{\textnormal{fast}}}
\newcommand{\rs}{\rh_{\textnormal{slow}}}
\newcommand{\rsplit}{\rh_{\textnormal{split}}}
\newcommand{\musplit}{\muh_{\textnormal{split}}}
\newcommand{\PP}[1]{\mathbb{P}\left\{{#1}\right\}} 
\newcommand{\Yset}{\mathcal{Y}}
\newcommand{\Ytrim}{\Yset_{\textnormal{trim}}}
\newcommand{\Ypred}{\Yset_{\textnormal{predict}}}
\newcommand{\Ypredfull}{\widetilde\Yset_{\textnormal{predict}}}
\newcommand{\Ysplit}{\Yset_{\textnormal{split}}}
\newcommand{\atrim}{\alpha_{\textnormal{trim}}}
\newcommand{\apred}{\alpha_{\textnormal{predict}}}
\newcommand{\eqnref}[1]{\eqref{eqn:#1}}
\newcommand{\betah}{\widehat{\beta}}
\newcommand{\norm}[1]{\lVert{#1}\rVert}
\newcommand{\bY}{\mathbf{Y}}
\newcommand{\bX}{\mathbf{X}}

\DeclareMathOperator*{\argmin}{arg\,min}

\begin{document}

\title{\Large\bfseries Trimmed Conformal Prediction for High-Dimensional Models}
\author[1]{Wenyu Chen}
\author[1]{Zhaokai Wang}
\author[1]{Wooseok Ha}
\author[1]{Rina Foygel Barber}
\affil[1]{Department of Statistics, University of Chicago}
\maketitle

\begin{abstract}
In regression, conformal prediction is a general methodology to construct prediction intervals in a distribution-free manner. Although conformal prediction guarantees strong statistical property for predictive inference, its inherent computational challenge has attracted the attention of researchers in the community. In this paper, we propose a new framework, called Trimmed Conformal Prediction (TCP), based on two stage procedure, a trimming step and a prediction step. The idea is to use a preliminary trimming step to substantially reduce
the range of possible values for the prediction interval, and then applying conformal prediction becomes far more efficient.
 As is the case of conformal prediction, TCP can be applied to any regression method, and further offers both statistical accuracy and computational gains. For a specific example, we also show how TCP can be implemented in the sparse regression setting. The experiments on both synthetic and real data validate the empirical performance of TCP.

\end{abstract}

\section{INTRODUCTION}

High-dimensional data are omnipresent in many scientific fields such as computational biology, and accordingly modern statistical inference has been evolving rapidly to develop tools that are valid in high-dimensional data. Though most research focuses explicitly on constructing and providing the confidence sets of the parameters of interest, they inevitably require strong assumptions on the data generation distribution or noise distribution. In contrast, recent work by \citet{Lei2016} studies a general framework for predictive inference under no distribution assumptions and only requiring the exchangeability of the training data and a new test data point. In their work, the main tool allowing the distribution free predictive inference is a method called ``conformal prediction'', first proposed by \citet{vovk2009,vovk1999,vovk2005} in the context of sequential classification and regression problems. To put it simply, conformal prediction is a way of exploiting the exchangeability of the residuals by treating them equally and sorting in the magnitudes---we provide a more detailed background in the next section. While conformal prediction is a powerful method to construct the prediction intervals in a distribution free manner, as pointed out by \citet{Lei2016}, it it usually computationally prohibitive. To avoid the computational issue, the authors also propose an alternative procedure based on sample splitting, called ``split conformal prediction''. Although the sample splitting method can reduce the computational cost drastically, the obvious drawback is losing the full power of the available sample, which is problematic
as it can substantially increase the width of the resulting prediction intervals.

In this paper, we address such issue and propose an efficient algorithm, Trimmed Conformal Prediction (TCP). TCP consists of a trimming step and a prediction step: in the trimming step, potential $y$ values for the test point are rapidly preprocessed to determine which $y$ values are most unlikely; these are then ``trimmed'' away. Then in the prediction step, we apply the conformal prediction algorithm only on this restricted set of $y$ values, to construct the prediction intervals more efficiently. The key idea of TCP is that by performing conformal prediction on a small set of candidate $y$ values, the computation can be improved significantly. Of course this comes at the price of losing some confidence in the prediction intervals, but in practice the loss is almost negligible. Our main theoretical result provides a simple and rigorous quantification of the loss of the coverage rate due to the trimming step.

To manifest the computational benefits of our algorithm, we also apply our method to the sparse linear model; here we use the lasso~\cite{tibshirani1996} as our prediction algorithm, and show that an easy trimming step can form a relatively small set of candidate $y$ values for the conformal prediction algorithm. In total, our algorithm requires far fewer calls of the lasso solver while losing almost no coverage rate in practice. \citet{hebiri2010} also considered conformal prediction based on the lasso estimator, but they do not take account for the change of the support when adding a new test point to the training data and thus their method does not offer a guarantee of finite sample coverage. On the other hand, our method guarantees the exact finite sample coverage while enjoying fast computation. We also demonstrate the superior performance of our method through a series of numerical experiments and real data applications.

\section{BACKGROUND}

\subsection{Sparse Models and the Lasso}

In the task of predicting the response $y$ at new feature point $x$,  linear regression is often preferred due to its simplicity and interpretability. In the high-dimensional regime where $p>n$, it is natural to assume that the regression coefficient is sparse. Under such assumption, the most popular approach is by solving the $\ell_1$-penalized least squares problem, also known as the lasso~\cite{tibshirani1996}, given by
\begin{equation}\label{eqn:lasso}
\betah = \argmin_{\beta\in\R^p}\left\{\frac{1}{2}\norm{\bY - \bX\beta}^2_2 + \lambda\norm{\beta}_1\right\},
\end{equation}
where $\bY\in\R^n$ is the response while  $\bX\in\R^{n\times p}$ contains the $p$ potential features.
One attractive property of the lasso~\eqnref{lasso} is that it selects only a few variables in the solution $\betah$, resulting a sparse model. While many results are known guaranteeing strong accuracy of the lasso
solution when the assumed sparse linear model is true or approximately true (e.g.~\cite{zhao2006,wainwright2009}), 
in practice we cannot know whether we are in such a setting, and assessing the predictive accuracy of the lasso without such assumptions is critical.

\subsection{Conformal Prediction}\label{sec:CPbackground}
Conformal prediction~\cite{vovk1999} seeks to provide reliable prediction intervals for 
new data points without assuming that the model used for prediction is necessarily
true. The only assumption is that the training and test data points are exchangeable.
Let $\widehat{\mu}:\mathcal{X}\rightarrow \R$ denote a 
{\em fitted regression function}, which is determined
by an unordered sample of $n+1$ many $(x,y)$ points, and maps
 $x$ values to predicted $y$ values.

Let $(X_1,Y_1),\dots,(X_n,Y_n)$ be our training data, and let $(X_{n+1},Y_{n+1})$
be the new test point where $X_{n+1}$ is known while $Y_{n+1}$ is the true but
unknown response value. 
The key idea is that if we were able to fit $\widehat{\mu}$ using this entire
sample of size $n+1$, $\{(X_i,Y_i):i=1,\dots,n+1\}$, then the exchangeability of the
$n+1$ data points implies that the residuals, $\rh_i = Y_i - \muh(X_i)$, are exchangeable
as well. In particular, the residual $\rh_{n+1}$ of the test point is equally
likely to rank anywhere in the list, and so the event that
$|\rh_{n+1}|$ is in the bottom $1-\alpha$ quantile of $|\rh_1|,\dots,|\rh_{n+1}|$
has probability at least $1-\alpha$. 

Of course, we cannot compute $\muh$ and $\rh$
since they depend on the unknown test point value $Y_{n+1}$. 
Instead, define
$\muh_y$ to be the fitted regression function using data points
$(X_1,Y_1),\dots,(X_n,Y_n),(X_{n+1},y)$ and let
\[\rh_y = (Y_1 - \muh_y(X_1),\dots,Y_n - \muh_y(X_n),y - \muh_y(X_{n+1}))\]
  be the residuals using
this function on this set of $n+1$ data points.
Note that if we happen to choose $y=Y_{n+1}$ we
recover the original values, i.e.~$\muh_{Y_{n+1}} = \muh$ and $\rh_{Y_{n+1}}=\rh$.
Then define
\[
\Yset = \Big\{y\in\R: \text{$|(\rh_y)_{n+1}|$ is in the 
bottom $1-\alpha$ quantile of $|(\rh_y)_1|,\dots,|(\rh_y)_{n+1}|$ }\Big\}.
\]
For $y=Y_{n+1}$, since $\rh_y = \rh$, we see that this statement must
hold with probability at least $1-\alpha$. That is, 
\[\PP{Y_{n+1}\in\Yset}\geq 1-\alpha,\]
giving the desired probability of coverage for the prediction interval $\Yset$.

\subsection{Split Conformal Prediction}
While the conformal prediction method described above gives the desired
statistical properties (namely, the correct coverage level without assuming any model
or distribution), computationally it can be quite challenging. The reason is that
the regression model must be fitted once for {\em each} value $y$ which is being 
considered---technically, we would have to fit a model $\muh_y$ for all $y\in\R$, although
in practice a large and fine grid of $y$ values may be used. If the model is computationally
demanding, this process can be inefficient. To address this, \citet{Lei2016} propose
a split conformal prediction procedure, where $\musplit$ is fitted on half the training sample,
$\{(X_i,Y_i):i\in I_1\}$ where $I_1\subset [n]$ is a set of size $|I_1| = n/2$, while the remaining half is used to obtain an empirical
distribution of the residuals. Specifically, defining residuals $(\rsplit)_i = Y_i - \musplit(X_i)$ for
$i\in I_2 = [n]\backslash I_1$, then the prediction interval is given by
\[
\Ysplit = \big\{y\in\R: \text{$|y - \musplit(X_{n+1})|$
 is in the bottom $(1-\alpha)$ quantile of $\{|(\rsplit)_i| : i\in I_2 \cup \{n+1\}\}$ }\big\}.
\]
As pointed out in~\cite{Lei2016} this set can simply be computed as the interval
\[\Ysplit = \musplit(X_{n+1}) \pm \rsplit^*\]
 where
$\rsplit^*$ is the $\lceil(1-\alpha)\cdot (n/2+1)\rceil$-st
 smallest value among $\{|(\rsplit)_i| : i\in I_2\}$.

This gives the same target coverage level, while only requiring one regression $\musplit$
to be computed; the drawback is that the intervals are likely to be wider, as the model
is fitted with only $n/2$ many samples and is therefore less accurate.

\section{METHOD}\label{sec:methods}
Before presenting our method, we motivate it with a simple example.
For the original form of conformal prediction, except in some special cases like linear
regression where the resulting interval can be computed in closed form, in general
we would need to choose some finite range $[-\ymax,\ymax]$ within which we select a grid
of $y$ values that we test for inclusion in $\Yset$. An intuitive choice for
$\ymax$ is to choose the largest value of the training data, $\ymax=\max_{1\leq i\leq n}|Y_i|$.
We can justify this theoretically by observing that, due to exchangeability,
\[\PP{|Y_{n+1}| > \max_{1\leq i\leq n} |Y_i|}\leq  \frac{1}{n+1},\]
and so 
\begin{equation}\label{eqn:cov_ymax}
\PP{Y_{n+1}\in\Yset\cap [-\ymax,\ymax]}\geq 1 - \alpha - \frac{1}{n+1},\end{equation}
giving nearly the nominal coverage level.
(Here $\Yset$ is the output of conformal prediction run on the full range $y\in\R$, as described
in Section~\ref{sec:CPbackground}.)

Our method, TCP, can be viewed as a generalization of this idea, and is summarized
with these two steps:
\begin{enumerate}
\item Trimming step: apply conformal prediction with a {\em fast but less accurate} 
method
to construct a  prediction interval $\Ytrim$ (which is wide, but generally will be 
much smaller than $[-\ymax,\ymax]$).
\item Prediction step: apply conformal prediction with a {\em slow but accurate}
regression model, working only over the restricted set $\Ytrim$.
\end{enumerate}
The goal is to obtain (nearly) the same accuracy as the second (slow) model,
while saving significant computation time by only fitting this slow regression model
over a restricted range $\Ytrim$.

To achieve this goal we consider two possibilities: trimming via a preliminary
conformal prediction step, or a split conformal prediction step.

\paragraph{Trimming with Conformal Prediction}
Suppose we are equipped with two regression algorithms, one 
fast (and inaccurate) and one slow (and highly accurate),
which each map a (unordered) data set of size $n+1$ to a 
predictive function, denoted as $\muf$ and $\mus$, respectively.
For any $y\in\R$, let $(\muf)_y$ denote the model $\muf$ fitted to data set
$(X_1,Y_1),\dots,(X_n,Y_n),(X_{n+1},y)$, and let $(\rf)_y$ denote the 
residuals of this model. Define $(\mus)_y,(\rs)_y$ analogously.

Our algorithm is given by the steps:
\begin{enumerate}
\item Trimming step:
\begin{multline*}
\Ytrim = \Big\{
y\in\R : \text{$|((\rf)_y)_{n+1}|$ is in the bottom $(1-\atrim)$ quantile}\\
\text{of $|((\rf)_y)_1|,\dots,|((\rf)_y)_{n+1}|$}\Big\}.
\end{multline*}
\item Prediction step:
\begin{multline*}
\Ypred = \Big\{
y\in\Ytrim : \text{$|((\rs)_y)_{n+1}|$ is in the bottom $(1-\apred)$ quantile}\\
\text{of $|((\rs)_y)_1|,\dots,|((\rs)_y)_{n+1}|$}\Big\}.
\end{multline*}
\end{enumerate}
Note the computational advantage: the slow model $(\mus)_y$
only needs to be fitted over $y\in\Ytrim$ (or rather, over a grid
of $y$ values covering $\Ytrim$), rather than over all $y\in\R$.
(In practice, the trimming step may need to 
be carried out over a finite grid of $y$ values, as well.)

In fact, we can consider a special case: if we choose the
fast method $\muf$ to simply predict a value of zero always, $\muf\equiv0$,
and set $\atrim = \frac{1}{n+1}$,
then this reduces to the method described earlier where conformal prediction
is computed only on the empirical range $[-\ymax,\ymax]$, since
the first step of our algorithm would 
compute $\Ytrim = [-\ymax,\ymax]$, and would then test only those
$y$ values in this range for inclusion into $\Ypred$.

\paragraph{Trimming with Split Conformal Prediction}
In some cases, we may prefer to use a split conformal prediction approach
for the trimming step. In this setting, consider again two regression
algorithms, $\musplit$ and $\mus$. We then carry out the steps:
\begin{enumerate}
\item Trimming step: apply split conformal prediction using $\musplit$, namely,
fit $\musplit$ to the first half of the training data,
then define
\[\Ytrim= \musplit(X_{n+1}) \pm \rsplit^*\]
where $\rsplit^*$ is the $\lceil (1-\atrim)(n/2+1)\rceil$-th smallest absolute residual
of $\musplit$ on the second half of the training data.
\item Prediction step: same as before.
\end{enumerate}
This procedure again offers a computational advantage
over simply applying the costly regression method $\mus$
 to a large grid of $y$ values (without a trimming step).
Of course, we also have the option of simply using split conformal prediction;
relative to this option, we instead have a {\em statistical} advantage---since
we are using the full sample size $n$ to build our predictive interval $\Ypred$,
we expect to have a narrower interval than if we had used the less accurate
model $\musplit$ which is fitted only on a sample of size $n/2$. We will see
these tradeoffs in practice in our experiments below.

\paragraph{Theoretical Guarantee}
As expected, our algorithm (in either form) offers a coverage guarantee that
combines the coverage levels of the two steps:
\begin{theorem}
If the data points $(X_1,Y_1)$, \dots, $(X_{n+1},Y_{n+1})$ are exchangeable,
then either version of the TCP method gives coverage level
\[\PP{Y_{n+1}\in\Ypred} \geq 1 - \atrim - \apred.\]
\end{theorem}
Before proving this result, we note that this result reduces to 
the coverage level given in~\eqnref{cov_ymax}
in the special case where we choose $\muf\equiv0$.
\begin{proof}
We see that $\Ytrim$ is simply the conformal prediction interval
determined by the fast regression model $\muf$ (in the first version of the method)
or the interval determined by split conformal inference with $\musplit$ (in the second version).
Furthermore, defining
\begin{multline*}
\Ypredfull = \Big\{
y\in\R : \text{$|((\rs)_y)_{n+1}|$ is in the bottom $(1-\apred)$ quantile}\\
\text{of $|((\rs)_y)_1|,\dots,|((\rs)_y)_{n+1}|$}\Big\},
\end{multline*}
this is the conformal prediction interval calculated with the slow model $\mus$.
Then, using existing results on conformal prediction (or split conformal prediction),
we know that $\PP{Y_{n+1}\in\Ytrim} \geq 1 - \atrim$ and $\PP{Y_{n+1} \in\Ypredfull}\geq 1 - \apred$.
By definition of our algorithm, we can also see that $\Ypred = \Ytrim \cap \Ypredfull$,
which proves the desired coverage level.
\end{proof}

\subsection{Application to Sparse Regression with the Lasso}
For the sparse linear model setting, as mentioned before, we would often
like to perform conformal prediction with the lasso~\eqnref{lasso} but cannot
afford to refit the lasso over a long list of $y$ values (e.g.~over a fine grid
spanning the range $[-\ymax,\ymax]$). 

We note here that the lasso is known to be easy to solve if the support and signs of 
the solution $\betah$ are known in advance---the solution then takes a closed form.
After fitting the lasso to some $y_*\in [-\ymax,\ymax]$, we might then hope that the support
and signs of $\betah$ would remain unchanged across this entire interval; this could easily be 
checked using the KKT optimality conditions for the lasso solution, and 
we would then avoid a second call to the lasso algorithm. 
Unfortunately, in practice, we have found that the support changes many times (even 
dozens of times) across this range even in simple simulations.

As an alternative, we run our two-stage algorithm with a more aggressive trimming step
to reduce the trial set farther.
 We now give the details.

\subsubsection{Trimming Step}
For the trimming step, 
 we consider two possibilities: 
\paragraph{Trimming via Ridge Regression} For the ridge regression option,
 the fast model $\muf$ is given by the penalized least squares regression
\[
\betah_{\textnormal{ridge}} = \argmin_{\beta\in\R^p}\left\{\frac{1}{2}\norm{\bY-\bX\beta}^2_2 + \frac{\rho}{2}\norm{\beta}^2_2\right\} \\= (\bX^\top \bX + \rho\mathbf{I}_p)^{-1}\bX^\top \bY,\]
where $\bX\in\R^{(n+1)\times p}$ and $\bY\in\R^{n+1}$ contain the full (training and test) data.
In this case we can solve for $\Ytrim$ in closed form, for the $n+1$ residuals can be written as $(\rf)_y = \mathbf{u}+\mathbf{v}y$, a linear function in $y$, where
\[
    \mathbf{u} =  \begin{bmatrix}
    \bY_{1:n}\\ 0\end{bmatrix}- \bX(\bX^\top \bX + \rho \mathbf{I}_p)^{-1}\bX^\top \begin{bmatrix}
    \bY_{1:n}\\ 0\end{bmatrix}, \quad
    \mathbf{v} =
    \begin{bmatrix}
\mathbf{0}_n\\  1\end{bmatrix} -
\bX(\bX^\top \bX + \rho \mathbf{I}_p)^{-1}X_{n+1}.
\]
For each $i$, then, the inequality $|((\rf)_y)_{n+1}|\leq|((\rf)_y)_i|$ holds in the interval enclosed by $\frac{\mathbf{u}_{n+1}-\mathbf{u}_i}{\mathbf{v}_i - \mathbf{v}_{n+1}}$ and $\frac{-\mathbf{u}_{n+1}-\mathbf{u}_i}{\mathbf{v}_i + \mathbf{v}_{n+1}}$.
Denote the smaller of these two values as $c_i$ and the larger as $d_i$; we can show that
$\max_i c_i \leq \min_i d_i$ always, since for all $i$,
the interval $[c_i,d_i]$ contains the value $y_* = -\mathbf{u}_{n+1}/\mathbf{v}_{n+1}$, which is the unique value for $y$ ensuring a zero residual value, i.e.~$((\rf)_y)_{n+1}=0$.
Finally, choose endpoints $c_{\textnormal{ridge}}\leq y_*\leq d_{\textnormal{ridge}}$ such
that $[c_i,d_i]\subseteq [c_{\textnormal{ridge}},d_{\textnormal{ridge}}]$
is violated for exactly $\atrim(n+1)$ many values $i\in[n]$; this gives the trimmed range, $\Ytrim = [c_{\textnormal{ridge}},d_{\textnormal{ridge}}]$. 
In particular, if $\atrim = \frac{1}{n+1}$, then we can simply set $\Ytrim = [\min_{i\in[n]} c_i,\max_{i\in[n]} d_i]$.

\paragraph{Trimming via Split Conformal Inference with Lasso}
Alternately, we can apply the lasso to half of the training data, indexed by $I_1\subset[n]$ with $|I_1|=n/2$, then use this
fitted model for trimming via split conformal inference. Specifically, define
\begin{equation}\label{eqn:LTSlasso}
\betah_{\textnormal{trim}} = \argmin_{\beta\in\R^p}
\left\{\frac{1}{2}\norm{\bY_{I_1}-\bX_{I_1}\beta}^2_2 + \lambda\norm{\beta}_1\right\}.\end{equation}
We then set 
\[\Ytrim = X_{n+1}^\top \betah_{\textnormal{trim}}  \pm \rsplit^*\]
where, as before, $\rsplit^*$ is the $\lceil (1-\atrim)(n/2+1)\rceil$-th smallest
value among
\[\{|Y_i - X_i^\top  \betah_{\textnormal{trim}} | : i \in I_2 = [n]\backslash I_1\}.\]
In particular, if we choose $\atrim = \frac{1}{n/2+1}$, then we obtain $\Ytrim = X_{n+1}^\top \betah_{\textnormal{trim}}  \pm\max_{i\in[n]} |Y_i - X_i^\top \betah_{\textnormal{trim}}|$.

\subsubsection{Prediction Step}
For the prediction step, we simply apply the conformal prediction algorithm using  lasso regression as the base method,
over only the restricted trial set $y\in\Ytrim$ (or, in practice, a fine grid of points over this set).
Specifically, writing
\[\betah_y = \argmin_{\beta\in\R^p} \left\{\frac{1}{2}\left\|\begin{bmatrix} \bY_{1:n}\\ y \end{bmatrix} - \bX\beta\right\|^2_2 + \lambda\norm{\beta}_1\right\},\]
we define
\begin{multline*}
\Ypred = \Big\{
y\in\Ytrim : \text{$|y - X_{n+1}^\top \betah_y|$ is in the bottom $(1-\apred)$ quantile of}\\
\text{$|Y_1 - X_1^\top\betah_y|,\dots,|Y_n - X_n^\top\betah_y|, |y - X_{n+1}^\top \betah_y|$}\Big\}.
\end{multline*}

\section{EMPIRICAL RESULTS}
We now test our method empirically,
both on simulated data and on real data from the Capital Bikeshare program in Washington DC. Code to reproduce our experiments is available online.\footnote{Available at \url{https://www.stat.uchicago.edu/\~rina/conformal.html}.}

\subsection{Simulation Studies}
We generate simulated data to compare the performance of four algorithms: 
conformal inference with the lasso on the trial set $[-\ymax,\ymax]$ (denoted as ``MaxTrim'' in our results); 
trimmed conformal inference using ridge regression to produce $\Ytrim$ (``RidgeTrim'');
trimmed conformal inference where $\Ytrim$ is obtained via split conformal inference with the lasso (``SplitTrim'');
and split conformal inference with the lasso~\cite{Lei2016} (``Split''). 
Note that MaxTrim, which performs conformal inference on the interval $[-\ymax,\ymax]$,
is essentially how conformal prediction is implemented in practice since, for methods such as lasso that must be rerun for each new $y$ value, we
cannot perform conformal prediction over the entire real line $y\in\R$.

In each case we aim for 90\% coverage. 
For the trimmed methods, we set $\atrim$ as small as possible, namely $\atrim = \frac{1}{n+1}$ (or $\atrim = \frac{1}{(n/2)+1}$ for trimming
with split conformal inference); see Section~\ref{sec:methods} for details on this. Since this is very small, 
we then set $\apred = 0.1$ to obtain (nearly) 90\% coverage.

\subsubsection{Settings}
We generate data from a linear model, $\bY = \bX \beta + \epsilon$, where:
\begin{itemize}
\item The $p$ features (the columns of $\bX$) are generated either with i.i.d.~$N(0,1)$ entries,
or are generated with high correlations, with the rows of $\bX$ drawn i.i.d.~from a $N(0,\Sigma)$ distribution,
where $\Sigma_{i,j}=0.9^{|i-j|}$;
\item The noise is generated either as $\epsilon_i\stackrel{\textnormal{i.i.d.}}{\sim} N(0,1)$ or $\epsilon_i\stackrel{\textnormal{i.i.d.}}{\sim} t(5)$, the $t$ distribution with 5
degrees of freedom;
\item The true signal $\beta$ is given $\beta_i = 2$ for $i\in S$, where the true support $S$ is a set of size $k$ chosen
at random, and $\beta_i = 0$ for $i\not\in S$.
\end{itemize}
Our methods will be run with desired level $\alpha = 0.1$.
We use $\lambda = \sqrt{n\log(p)}$ for normal noise case, and $\lambda = \sqrt{\frac{5}{3} n\log(p)}$ for $t$ noise case (when we use split conformal inference, and the effective sample size is $n/2$, we use $\lambda = \sqrt{(n/2)\log(p)}$ or $\sqrt{\frac{5}{3} (n/2)\log(p)}$).

\subsubsection{Results}
Table~\ref{tab:sim_results} displays our results for the setting $n=200$, $p=2000$, $k=10$, across both feature and noise models;
we record the coverage level for each method, as well as the width of the trial set $\Ytrim$ (which should be interpreted as the computational cost---note
that this is not relevant for split conformal inference, which only ever fits the model once), and the width of the resulting prediction interval $\Ypred$.
Results are averaged over 500 trials.

\begin{figure}[!t]
\centering
\includegraphics[width=3in]{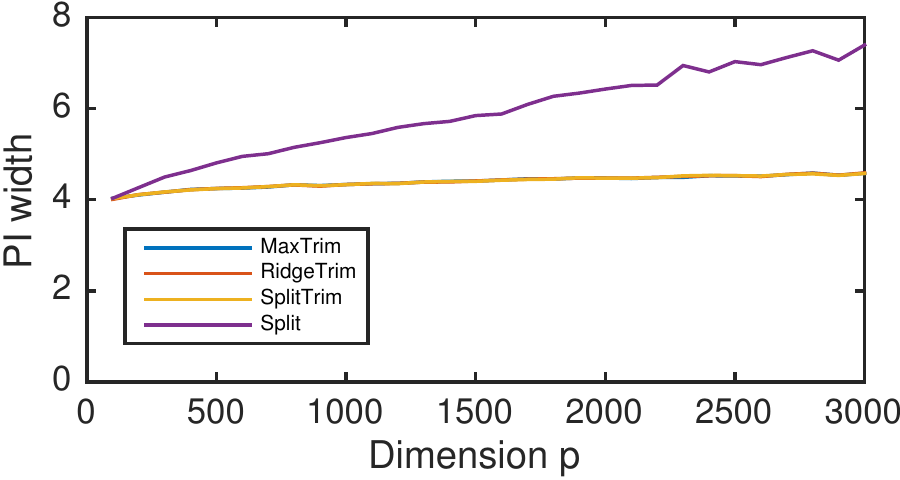}
\caption{Width of prediction intervals (PIs) for all methods, as the number of features $p$ varies. The results for MaxTrim, RidgeTrim, and SplitTrim
are nearly identical and cannot be distinguished in the plot.}
\label{fig:vary_p_PIwidth}
\end{figure}

\begin{figure}[!t]
\centering
\includegraphics[width=3in,clip=TRUE]{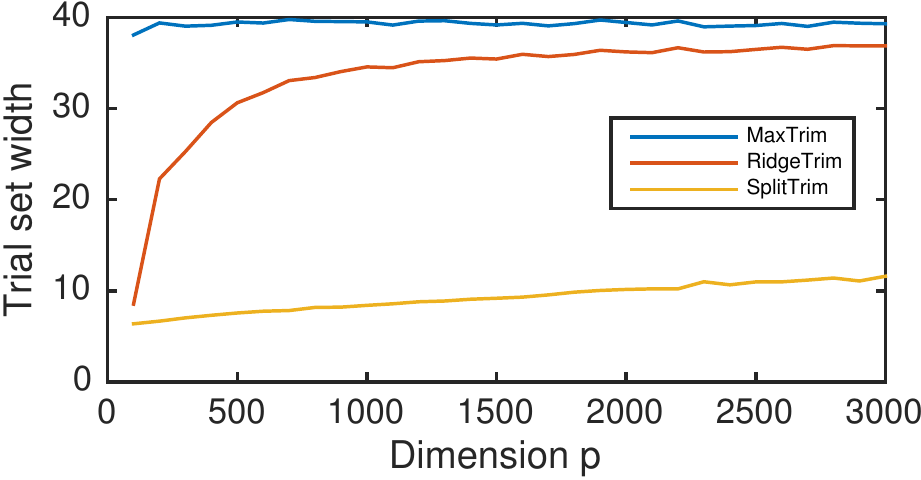}
\caption{Width of trial sets for all methods, as the number of features $p$ varies; this width is a proxy for computational cost. (The split conformal inference method
is not included as it only fits the lasso model once in any setting.)}
\label{fig:vary_p_trialwidth}
\end{figure}

\begin{figure}[!t]
\centering
\includegraphics[width=3in,clip=TRUE]{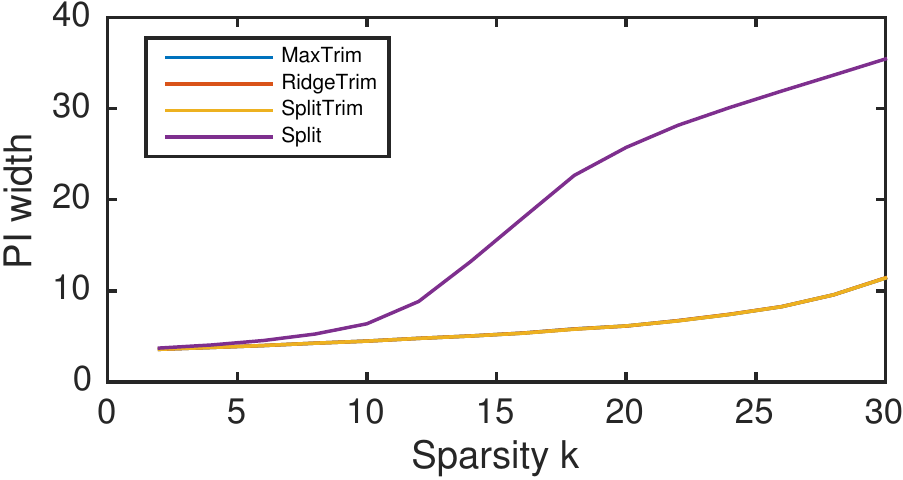}
\caption{Width of prediction intervals (PIs) for all methods, as the size $k$ of the true model varies. The results for MaxTrim, RidgeTrim, and SplitTrim
are nearly identical and cannot be distinguished in the plot.}
\label{fig:vary_k_PIwidth}
\end{figure}

\begin{figure}[!t]
\centering
\includegraphics[width=3in,clip=TRUE]{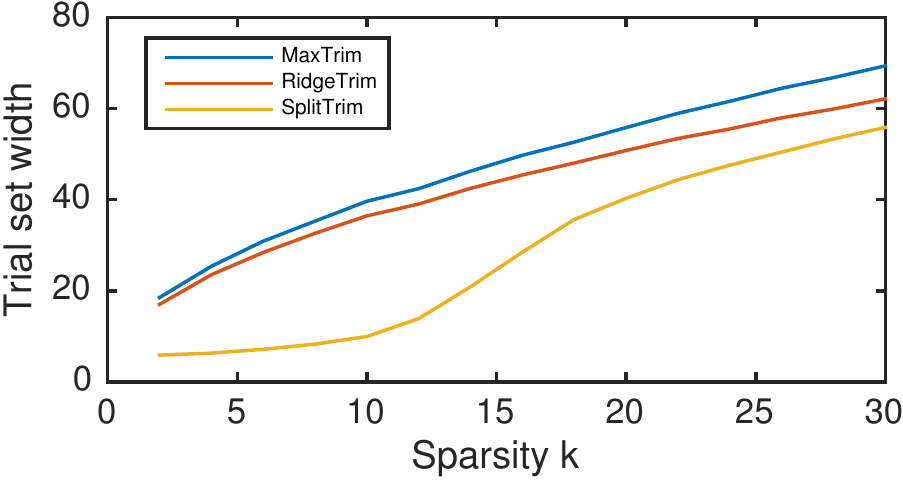}
\caption{Width of trial sets for all methods, as the size $k$ of the true model varies; this width is a proxy for computational cost. (The split conformal inference method
is not included as it only fits the lasso model once in any setting.)}
\label{fig:vary_k_trialwidth}
\end{figure}

In Figures~\ref{fig:vary_p_PIwidth} and~\ref{fig:vary_p_trialwidth}, we show the widths of $\Ypred$ and $\Ytrim$ when the number of features
varies, $p=100,200,\dots,3000$, with $n=200$ and $k=10$ fixed. We do the same in Figures~\ref{fig:vary_k_PIwidth} and~\ref{fig:vary_k_trialwidth}
when instead $n=200, p=2000$ are fixed and we vary the true model size, $k=2,4,\dots,30$. For these figures, we show results only for the normal noise model and
the uncorrelated features model.
Results are averaged over 500 trials.

{
\begin{table}[!t]\footnotesize
\centering
\begin{tabular}{|@{\,}c@{\,}|@{\,}c@{\,}|c|c|c|c|}
\hline
&&&PI&Trial set&Coverage\\
&&&width&width&(\%)\\\hline
\multirow{8}{*}{\begin{sideways}Uncorr.~features\end{sideways}}
&\multirow{4}{*}{\begin{sideways}$N(0,1)$\end{sideways}}
&MaxTrim & 4.39 &  38.06 &   90.4 \\
&&RidgeTrim &4.39 &  36.22 &   90.4 \\
&&SplitTrim &4.39 &  10.06 &   90.2 \\
&&Split & 6.35 & ---  &  89.6 \\\cline{2-6}
&\multirow{4}{*}{\begin{sideways}$t(5)$\end{sideways}}
&MaxTrim & 5.54 & 38.51 & 92.0\\
&&RidgeTrim & 5.54&36.99 & 92.0 \\
&&SplitTrim & 5.54 & 13.76 & 92.0 \\
&&Split & 8.27 &  ---& 91.8\\\hline
\multirow{8}{*}{\begin{sideways}High corr.~features\end{sideways}}
&\multirow{4}{*}{\begin{sideways}$N(0,1)$\end{sideways}}
&MaxTrim &4.41 & 38.18 & 91.2 \\
&&RidgeTrim & 4.40 & 26.50  &  91.0 \\
&&SplitTrim & 4.41 &  9.09   &  91.0 \\
&&Split &  5.77 & --- & 89.6\\\cline{2-6}
&\multirow{4}{*}{\begin{sideways}$t(5)$\end{sideways}}
&MaxTrim & 5.55& 39.02& 90.0 \\
&&RidgeTrim & 5.54 & 27.88 & 89.8 \\
&&SplitTrim & 5.54 & 12.53 & 90.0 \\
&&Split & 7.34&  ---&91.8\\\hline
\end{tabular}
\caption{Results of the simulated data experiment, for uncorrelated or highly correlated
features, and for $N(0,1)$ or $t(5)$ distributed noise, for settings $n=200,p=2000,k=10$.
Desired coverage level is $1-\alpha = 90\%$.
``Trial set width'' is the width
of $\Ytrim$ for the trimmed conformal prediction methods, and should be viewed as a proxy for computational cost.}
\label{tab:sim_results}
\end{table}
}

From the simulation, we see that the three TCP methods (MaxTrim, RidgeTrim, and SplitTrim) give substantially
narrower prediction intervals than split conformal prediction. With more features $p$ or a larger true model size $k$, 
the gap between split conformal inference and the TCP methods increases---this is intuitive, as using the full sample size 
becomes more important for accuracy as the problem becomes more difficult.

 In terms of computational cost, among the three trimmed methods, SplitTrim is always
 the most aggressive, giving the sharpest (smallest) $\Ytrim$. RidgeTrim is more conservative for this data,
 but nonetheless improves over MaxTrim.
 
 The coverage is near the desired level of 90\% for all methods in all settings.
  
\subsection{Real Data Analysis}
In this data analysis, we study the rental records of the Capital Bikeshare program in Washington DC.\footnote{Capital
Bikeshare data is publicly available at \url{https://www.capitalbikeshare.com/trip-history-data}.} The data consists of bike
rental records noting the date, time, and locations of the rental and return of each bicycle used in the program.
We aggregate this data to record only the total number of rentals originating at each station on each day.
We will choose one station to use as the response and the rest will be features; that is, we will predict the number of rentals originating
at one station on a given day, as a function of the number of rentals at each of the other stations.

As the business expands, more bike stations are added. Therefore, we only pick a period of time when no new stations appears:
the 93 days from Nov.~7,~2010 through Feb.~7,~2011. There are 107 stations with activity during these days.
Choosing one station $j_*\in \{1,\dots,107\}$ as the response and one day $i_*\in\{1,\dots,93\}$,
 as the test point, we therefore have $p=106$ many stations given by the remaining features $\{1,\dots,107\}\backslash\{j_*\}$,
and $n=92$ training data points given by the remaining days $\{1,\dots,93\}\backslash\{i_*\}$. 

We set $\alpha = 0.1$, and repeat our experiment with each station as the response, i.e.~$j_* = 1, 2,\dots, 107$, and average our results over these 107 runs.
For choosing $i_*$, we  consider two settings:\begin{itemize}
\item ``Random day'' setting. We pick a random day $i_*$ from the 93 data points. Results are averaged over 10 randomly
selected days.
\item ``Last day'' setting. We pick the most recent day as the test point, i.e.~$i_*=93$. This is more practical for real time analysis and prediction,
since in practice
we would typically use past data (from the first 92 days) to predict upcoming events (the present day, i.e.~the 93rd day). 
However,
as behavior patterns may change over time, this setting violates the exchangeability assumption, and the coverage
properties of the conformal inference methods may not hold.
\end{itemize}

Table~\ref{tab:realdata} shows the average coverage rates for the four algorithms under the two different test methods.

\begin{table}[h]
\centering
\begin{tabular}{|@{\,}c@{\,}|c|c|c|c|}
\hline
&&PI&Trial set&Coverage\\
&&width&width&(\%)\\\hline
\multirow{4}{*}{\begin{sideways}\scriptsize Random day\end{sideways}}
&MaxTrim & 11.09& 37.21& 92.1\\
&RidgeTrim & 11.08& 21.40& 92.1\\
&SplitTrim & 10.96& 18.03& 92.1 \\
&Split & 11.47& ---& 88.3\\\hline
\multirow{4}{*}{\begin{sideways}\scriptsize Last day\end{sideways}}
&MaxTrim & 11.36& 37.19& 86.0\\
&RidgeTrim & 11.38& 22.34& 86.0\\
&SplitTrim& 11.14& 18.15 & 86.0\\
&Split & 11.43& ---& 85.0\\\hline
\end{tabular}
\caption{Results for the Capital Bikeshare data experiment, for the ``random day'' and the ``last day'' settings.
Desired coverage level is $1-\alpha = 90\%$.
}
\label{tab:realdata}
\end{table} 



Despite slight numerical differences, in the ``random day'' setting, the
average coverage rates are all close to the desired level $1-\alpha = 90\%$.
For the ``last day'' setting, however, the coverage level is noticeably lower,
presumably due to the violation of exchangeability in this setting.

In this experiment, the prediction intervals are all quite similar in length
(however, it is worth noting that the split conformal inference method
achieves slightly lower coverage at roughly the same interval length compared
to the other methods), while the trimmed conformal prediction methods indeed
show substantial reduction in computation time (as measured by the length of the trial 
set, $\Ytrim$) relative to conformal prediction applied to the broad interval $[-\ymax,\ymax]$.

\section{DISCUSSION}
In this paper, we present a fast two-stage algorithm for conformal prediction, called Trimmed Conformal Prediction (TCP). In the trimming step, the most unlikely $y$ values are trimmed away quickly, and in the prediction step, the conformal prediction algorithm is applied over this reduced range of potential $y$ values. Our empirical results on simulated and real data show that TCP achieves computational gains over conformal prediction without the trimming step, while offering sharper prediction intervals (a statistical advantage) compared to split conformal prediction. This is highly desirable in the high-dimensional data analysis, where we are faced with both statistical and computational challenges.


\bibliographystyle{plainnat}
\bibliography{bib}

\end{document}